\documentclass[12pt]{amsart}
\usepackage[utf8]{inputenc}
\usepackage[T1]{fontenc}
\usepackage[all]{xy}
\usepackage{lipsum}
\usepackage{url}
\usepackage{tikz}
\usepackage{stackrel}
\usepackage{color}

\usepackage{amsmath,amsthm,amssymb}
\newcommand{\aspas}[1]{``{#1}''}
\usepackage{xcolor}
\usepackage{graphicx}

\newtheorem{theorem}{Theorem}[section]
\newtheorem{lemma}[theorem]{Lemma}
\newtheorem{example}[theorem]{Example}
\newtheorem{proposition}[theorem]{Proposition}
\theoremstyle{definition}
\newtheorem{definition}[theorem]{Definition}

\newtheorem{remark}[theorem]{Remark}
\newtheorem{corollary}[theorem]{Corollary}

\numberwithin{equation}{section}

\begin{document}

\vspace{0.5in}

\renewcommand{\bf}{\bfseries}
\renewcommand{\sc}{\scshape}
\vspace{0.5in}

\title[Configuration sets for groups]%
{Configuration sets for groups}

\author[C. A. Ipanaque Zapata]{Cesar A. Ipanaque Zapata}
\address[C. A. Ipanaque Zapata]{Departamento de Matem\'atica, IME-Universidade de S\~ao Paulo, Rua do Mat\~ao, 1010 CEP: 05508-090, S\~ao Paulo, SP, Brazil}
\email{cesarzapata@usp.br}

\author[Alex Campos]{Alex Freitas de Campos}
\address[Alex Freitas de Campos]{ICMC -- Universidade de São Paulo, Av. Trab. São Carlense, 400, São Carlos, SP, Brazil}
\email{afcampos@icmc.usp.br}


\thanks {The first author would like to thank grant\#2022/16695-7 and grant\#2016/18714-8, S\~{a}o Paulo Research Foundation (FAPESP) for financial support.}

\begin{abstract} We develop practical tools for analyzing the configuration set $F(G,k)$ of $k\geq 2$ distinct elements in a group $G$. We apply our results to design homogeneous linear systems in $\mathbb{F}_q$ that admit nontrivial solutions. Furthermore, we study the connectivity of Cayley graphs of the form $\mathrm{Cay}(G^k,F(G,k))$. In addition, we consider the configuration set of $k\geq 2$ distinct non-identity elements in $G$. 
\end{abstract}

\maketitle


\section{Introduction}
In this article, the term \aspas{homomorphism} refers to a group homomorphism. 

\medskip Let $G$ be an abstract group (without additional  structure), and let $k\geq 2$ be an integer. The $k$-th direct product of $G$, denoted by $G^k=G\times\cdots\times G$ ($k$ times), is equipped with the component-wise operation: \[(g_1,\ldots,g_k)\cdot(h_1,\ldots,h_k)=(g_1h_1,\ldots,g_kh_k).\] Notation: $$
\begin{array}{rcl}
g_1 & \cdots & g_k\\
h_1 & \cdots & h_k\\ \hline 
g_1h_1&\cdots & g_kh_k \\
\end{array}
$$ The ordered \textit{configuration set} of $G$ is the subset of $G^k$ defined by \[F(G,k)=\{(g_1,\ldots,g_k)\in G^k| g_i\neq g_j \text{ for all } i\neq j\}.\] More generally, one can consider the configuration set $F(X,k)$ of any set $X$. By convention $F(X,1)=X$. 

\medskip In algebraic topology, the set $F(G,k)$ has been extensively studied  when $G$ is a topological space $-$ commonly referred to as the \textit{ordered configuration space}~\cite{fadell1962configuration} $-$ and particularly when $G$ is a manifold~\cite{idrissi2022}. However, when $G$ is treated as an abstract group endowed with the discrete topology, $F(G,k)$ becomes a discrete space with limited topological interest, and its algebraic or combinatorial structure has received little to no attention. This work takes a first step in that direction.    

\medskip Note that $F(G,k)\neq\emptyset$ if and only if the order of the group satisfies $| G| \geq k$. Throughout this work, we assume that $| G| \geq k\geq 2$. Observe that the number of elements in $F(G,k)$ is \begin{align}\label{order-conf}
    |G|(|G|-1)(|G|-2)\cdots(|G|-k+1)&=|G|!/(|G|-k)!
\end{align} whenever $k\leq |G|<\infty$.

\medskip It is important to note that $F(G,k)$ is not a subgroup of the direct product $G^k$ (it does not even contain the identity element $(1,\ldots,1)$, where $1$ denotes the identity of $G$). 

\medskip The intersection of a family of subgroups of $G$ is again a subgroup. Therefore, for any subset $X\subseteq G$, there exists a smallest subgroup of $G$ containing $X$, called the \textit{subgroup generated} by $X$ and denoted by $\langle X\rangle$. Recall that \[\langle X\rangle=\{x_1\cdots x_n:~x_i\in X\cup X^{-1} \text{ for } 1\leq i\leq n, ~n\geq 1 \},\] where $X^{-1}=\{x^{-1}:~x\in X\}$. If $G=\langle X\rangle$, then $X$ is called a \textit{generating set} of $G$. 

\begin{remark}\label{invariance-conf}
If $f:G\to H$ is an epimorphism and $X$ is a generating set of $G$, then the image $f(X)$ is a generating set of $H$. In particular, if $f:G\to H$ is an isomorphism, then the induced map \[f^k:G^k\to H^k,~(g_1,\ldots,g_k)\mapsto (f(g_1),\ldots,f(g_k))\] is also an isomorphism, and satisfies $f^k(F(G,k))=F(H,k)$. Therefore, $F(G,k)$ is a generating set of $G^k$ if and only if $F(H,k)$ is a generating set of $H^k$.
\end{remark}

\medskip Given a group $G$ (finite or infinite) with  $\infty\geq |G|\geq k\geq 2$, we address the following natural questions: Is the configuration set $F(G,k)$ a generating set for $G^k$? If not, what subgroup is generated by $F(G,k)$? To the best of our knowledge, these questions have not been addressed in the existing literature. In this paper, we begin the study of these questions.  

\medskip The main results of this paper are as follows: 
\begin{itemize}
    \item We study whether the configuration set $F(G,k)$ is a generating set of $G^k$ under certain conditions (Theorem~\ref{the-conf-set-thm}).
    \item For odd $n\geq 3$, we show that the configuration set $F(D_n,2n)$ of the dihedral group $D_n$ of order $2n$ is a generating set of $D_n^{2n}$ (Theorem~\ref{ndihedral-group}).
    \item We provide an algebraic-combinatorial characterization of $F(G,k)$ (Theorem~\ref{charact-k2}).
    \item We develop a method for constructing homogeneous linear systems over  $\mathbb{F}_q$ that admit nontrivial solutions (Proposition~\ref{nontrivial-solution}). 
    \item We study the connectivity of the Cayley graphs $\mathrm{Cay}(G^k,F(G,k))$ (Theorem~\ref{cayley-conf}).
    \item We study the configuration set $F(G\setminus\{0\},k)$ for an abelian group $G$ (Theorem~\ref{thm:punctured-group}). 
\end{itemize}

This paper is organized into two sections. In Section~\ref{conf-set}, we study the set $F(G,k)$ and introduce new concepts that support this notion. For instance, we define the norm of $G^k$ (Definition~\ref{def:norm}). Section~\ref{sec:applications} presents new insights into group theory and combinatorics. For example, we introduce the configuration property (Definition~\ref{defn:cp}) and study the Cayley graphs of the form $\mathrm{Cay}(G^k,F(G,k))$. We also examine the configuration set $F(G\setminus\{0\},k)$. In Remark~\ref{rem:future-work}, we propose directions for future work.

\section{The configuration set $F(G,k)$}\label{conf-set}
In this section, we study whether the configuration set $F(G,k)$ is a generating set of $G^k$ (see Theorem~\ref{the-conf-set-thm}). A key concept in showing that $F(G,k)$ is not a generating set of $G^k$ is the notion of a norm introduced in Definition~\ref{def:norm}. 

\medskip Note that $F(G,k)$ is \textit{symmetric}, i.e., $F(G,k)=F(G,k)^{-1}$. Let $S_k$ denote the symmetric group on $k$ elements. For $\sigma\in S_k$ and $(g_1,\ldots,g_k)\in G^k$, we write \[\sigma(g_1,\ldots,g_k)=(g_{\sigma(1)},\ldots,g_{\sigma(k)}).\] In addition, note that if $(g_1,\ldots,g_k)\in F(G,k)$, then   $\sigma(g_1,\ldots,g_k)\in F(G,k)$ for any $\sigma\in S_k$. 

\medskip We now present the following example.

\begin{example}
\noindent\begin{enumerate}
    \item[(1)] Let $G=\mathbb{Z}_2=\{\overline{0},\overline{1}\}$ be the integers modulo 2, and let $k=2$. Then, \[F(\mathbb{Z}_2,2)=\{(\overline{0},\overline{1}),(\overline{1},\overline{0})\}\subseteq \mathbb{Z}_2^2.\] Furthermore, we have $(\overline{0},\overline{0})\notin F(\mathbb{Z}_2,2)$ and $(\overline{0},\overline{1})+(\overline{1},\overline{0})=(\overline{1},\overline{1})\notin F(\mathbb{Z}_2,2)$.
    \item[(2)] Let $G=\mathbb{Z}$ be the group of integers, and let $k=2$. Then, \[F(\mathbb{Z},2)=\{(m,n)\in\mathbb{Z}\times \mathbb{Z}:m\neq n\}\subseteq \mathbb{Z}^2.\] Again, we have $(0,0)\notin F(\mathbb{Z},2)$ and $(0,1)+(1,0)=(1,1)\notin F(\mathbb{Z},2)$.
\end{enumerate}
\end{example}

A key concept in showing that $F(G,k)$ is not a generating set of $G^k$ is the notion of a norm, defined as follows.

\begin{definition}[Norm]\label{def:norm} Let $G$ be a group and let $k\geq 2$ be an integer. The \textit{norm} on $G^k$ is the map $\|\cdot\|\colon G^k\to G$ defined by \[\| g\|=g_1\cdots g_k\] for any $g=(g_1,\ldots,g_k)\in G^k$. Note that the norm $\|\cdot\|$ is surjective. Indeed, it admits a \textit{global cross-section}, that is, a homomorphism $s:G\to G^k$ such that $\|\cdot\|\circ s=\mathrm{id}_G$. For example, $s(g)=(g,1,\ldots,1)$ for any $g\in G$. 
\end{definition}

Note that $G$ is abelian if and only if the norm $\|\cdot\|\colon G^k\to G$ is a group homomorphism, i.e., \[\| gh\|=\| g\|\, \| h\| \] for all $g,h\in G^k$ (for some, and thus for any, $k\geq 2$).

\medskip A general connection between the norm and configurations is the following statement.

\begin{proposition}\label{prop:ker-conf-incl}
  Let $G$ be a group and $k\geq 3$ such that $|G|\geq k$. We have \[\|\cdot\|^{-1}(1)\subseteq \langle F(G,k)\rangle.\]  
\end{proposition}
\begin{proof}
    Suppose $g=(g_1,\ldots,g_k)\in G^k$ such that $g_1g_2\cdots g_k=1$, i.e, $g_1=g_k^{-1}\cdots g_2^{-1}$. Hence, $$
\begin{array}{rccccccl}
g_k^{-1} & 1 & 1 & 1 & \ldots & 1 & 1 & g_k\\
g_{k-1}^{-1}& 1 & 1 & 1 & \ldots & 1 & g_{k-1} & 1\\
\vdots & \vdots & \vdots & \vdots & \ddots & \vdots & \vdots & \vdots \\
g_3^{-1} & 1 & g_{3} & 1 & \ldots & 1 & 1 & 1\\ 
g_2^{-1} & g_2 & 1 & 1 & \ldots & 1 & 1 & 1\\ \hline 
g_1 & g_2 & g_3 & g_4 & \ldots & g_{k-2} & g_{k-1} & g_k\\
\end{array} 
$$  
We claim that $(g^{-1},g,1,\ldots,1)\in \langle F(G,k)\rangle$ for any $g\in G$. It shows that $\|\cdot\|^{-1}(1)\subseteq \langle F(G,k)\rangle$. To prove the claim, if $g=1$, it is trivially true. If $g\neq 1$, choose elements $h_1,\ldots,h_{k-3}\in G\setminus\{g,g^{-1},1\}$ with $h_i\neq h_j$ for $i\neq j$ (this is possible since $| G| \geq k\geq 3$). Then, \[(g^{-1},g,1,\ldots,1)=(g^{-1},1,g,h_1,\ldots,h_{k-3})\cdot (1,g,g^{-1},h_1^{-1},\ldots,h_{k-3}^{-1}),\] and both factors lie in $F(G,k)$. Therefore, $(g^{-1},g,1,\ldots,1)\in \langle F(G,k)\rangle$.
\end{proof}

\medskip On the other hand, for $k\geq 2$, note that the subset \[E_k=\{(g_1,1,\ldots,1),(1,g_2,1,\ldots,1),\ldots,(1,\ldots,1,g_k):~g_i\in G\}\] is a generating set of $G^k$ for any group $G$. For each $g\in G$, define \[e_j(g)=(1,\ldots,1,g,1,\ldots,1)\in G^k\] where $g$ appears in the $j$-th coordinate and all other coordinates are $1$. Recall that, given $\sigma\in S_k$ and $(g_1,\ldots,g_k)\in G^k$, we write $\sigma(g_1,\ldots,g_k)=(g_{\sigma(1)},\ldots,g_{\sigma(k)})$. If $(g_1,\ldots,g_k)\in F(G,k)$, then $\sigma(g_1,\ldots,g_k)\in F(G,k)$ for any $\sigma\in S_k$. We now prove the following result.

\begin{lemma}\label{lem:ejg} Let $G$ be any group and $k\geq 2$ be an integer. 
\begin{enumerate}
    \item[(1)] Let $g\in G$. If $e_j(g)\in \langle F(G,k)\rangle$ for some $j\in\{1,\ldots,k\}$, then $e_i(g)\in \langle F(G,k)\rangle$ for all $i\in\{1,\ldots,k\}$. 
     \item[(2)] Suppose $G=\langle X\rangle$ is a group generated by a set $X$. If there exists $j\in\{1,\ldots,k\}$ such that $e_j(x)\in \langle F(G,k)\rangle$ for every $x\in X$, then $E_k\subseteq \langle F(G,k)\rangle$. In particular, the configuration set $F(G,k)$ is a generating set of $G^k$.
\end{enumerate} 
\end{lemma}
\begin{proof}
\noindent\begin{enumerate}
    \item[(1)] The case $i=j$ is trivial. Assume $i\neq j$. Let $\sigma_{i,j}$ be the transposition that swaps the $i$-th and $j$-th coordinates and leaves the others fixed. Suppose $e_j(g)=c_1\cdots c_m$, where each $c_\ell\in  F(G,k)$. Then, \begin{eqnarray*}
        e_i(g)&=&\sigma_{i,j}e_j(g)\\
        &=&(\sigma_{i,j}c_1)\cdots(\sigma_{i,j}c_m),
    \end{eqnarray*} and each $(\sigma_{i,j}c_\ell)\in F(G,k)$. Therefore, $e_i(g)\in \langle F(G,k)\rangle$.
    \item[(2)] Let $g\in G$, and suppose $g=x_1\cdots x_m$ with each $x_1,\ldots,x_m\in X\cup X^{-1}$. Note that if $x\in X$, then $(x^{-1},1,\ldots,1)=(x,1,\ldots,1)^{-1}\in \langle F(G,k)\rangle$, since $(x,1,\ldots,1)\in \langle F(G,k)\rangle$ by part (1). Therefore, \[(g,1,\ldots,1)=(x_1,1,\ldots,1)\cdots (x_m,1,\ldots,1)\] with each factor  in $\langle F(G,k)\rangle$. Hence, $(g,1,\ldots,1)\in \langle F(G,k)\rangle$. Then, by part (1), we get  $e_i(g)\in \langle F(G,k)\rangle$ for all $i\in\{1,\ldots,k\}$ and all $g\in G$, so $E_k\subseteq \langle F(G,k)\rangle$.  
\end{enumerate}  
\end{proof}

 The first main result of this section is the following.

\begin{theorem}\label{the-conf-set-thm}
Let $G$ be a group. 
\begin{enumerate}
    \item[(1)] If $k=2$ or $| G| \geq k+1\geq 4$, then $F(G,k)$ is a generating set of $G^k$. In particular, if $G$ is infinite, then the subgroup generated by  $F(G,k)$ is equal to $G^k$ for any $k\geq 2$.
    \item[(2)] If $G$ is a finite abelian group with $| G| =k\geq 3$, then $F(G,k)$ is not a generating set of $G^k$. Furthermore, \begin{align*}
    \langle F(G,k)\rangle&= \|\cdot\|^{-1}(\langle a\rangle),\quad \text{for any $k\geq 3$,}
\end{align*} where $a=\sum_{g\in G}g$ is an element of order at most 2. For instance, $\langle F(G,k)\rangle= 
   \mathrm{Ker}(\|\cdot\|)$ whenever $k$ is odd.
\end{enumerate}
\end{theorem}
\begin{proof}
    \noindent\begin{enumerate}
    \item[(1)] The case $k=2$ follows from Lemma~\ref{lem:ejg}(2), using the fact that $(g,1)\in F(G,2)$ for any $g\in G\setminus\{1\}$. 
    
    Now, suppose $k\geq 3$. By Lemma~\ref{lem:ejg}(2), it suffices to check that $(g,1,1,\ldots,1)\in \langle F(G,k)\rangle$ for every $g\in G$. If $g=1$, it is trivially true. If $g\neq 1$, choose elements  $h_1,\ldots,h_{k-1}\in G\setminus\{g,1\}$ with $h_i\neq h_j$ for $i\neq j$ (this is possible since $| G| \geq k+1$ and $k\geq 3$). Then, \[(g,1,1,\ldots,1)=(g,h_1,h_2,\ldots,h_{k-1})\cdot (1,h_1^{-1},h_2^{-1},\ldots,h_{k-1}^{-1}),\] and both factors lie in $F(G,k)$. Therefore, $(g,1,1,\ldots,1)\in \langle F(G,k)\rangle$.
        \item[(2)] Suppose $G=\{g_1,g_2,\ldots,g_{k}\}$ with  $g_1=0$ as the identity. Note that the configuration set $$F(G,k)=\{(g_{\sigma(1)},\ldots,g_{\sigma(k)}):~\sigma\in S_k\},$$ so any $g\in F(G,k)$ satisfies $\| g\|=g_2+\cdots+g_k$. Then, any element in the subgroup $\langle F(G,k)\rangle$ has norm in $\langle g_2+\cdots+g_k\rangle$, that is, \[\langle F(G,k)\rangle\subseteq  \|\cdot\|^{-1}(\langle g_2+\cdots+g_k\rangle).\] Observe that $(g_2+\cdots+g_k)+(g_2+\cdots+g_k)=0$, since $\{g_2,\ldots,g_k\}=\{-g_2,\ldots,-g_k\}$. Thus, $a:=g_2+\cdots+g_k$ has order at most 2 (for instance it has order 1 whenever $k$ is odd). Then $\langle a\rangle$ is a proper subgroup of $G$ because $|G|=k\geq 3$. Hence, $F(G,k)$ is not a generating set of $G^k$.

On the other hand, by Proposition~\ref{prop:ker-conf-incl}, we obtain \[\mathrm{Ker}(\|\cdot\|)\subseteq \langle F(G,k)\rangle\] for any $k\geq 3$. Furthermore, let $c=(c_1,\ldots,c_k)\in G^k$ such that $c_1+\cdots+c_k=a=g_1+\cdots+g_k$. Then \[c-g\in \mathrm{Ker}(\|\cdot\|)\subseteq \langle F(G,k)\rangle\] where $g=(g_1,\ldots,g_k)\in F(G,k)$. Hence, $c\in \langle F(G,k)\rangle$. It implies that \[\|\cdot\|^{-1}(\langle a\rangle)\subseteq \langle F(G,k)\rangle\] for any $k\geq 3$. Therefore, \[\langle F(G,k)\rangle=\|\cdot\|^{-1}(\langle a\rangle)\] for any $k\geq 3$.  

For instance, \begin{align*}
    \langle F(G,k)\rangle&= 
   \mathrm{Ker}(\|\cdot\|)  
\end{align*} whenever $k$ is odd.
    \end{enumerate}
\end{proof}

\begin{remark}
Let $G$ be a finite abelian group with $|G| =k\geq 3$. For $k$ even, the sum \[a:=\sum_{g\in G}g\] could be zero (see Example~\ref{exam:F4}).
\end{remark}

As a direct consequence of Theorem~\ref{the-conf-set-thm}(2),  we have the following result.

\begin{corollary}
Let $G$ be a finite group with $k=|G|\geq 3$. 
\begin{enumerate}
    \item[(1)] If $F(G,k)$ is a generating set of $G^k$, then $G$ is non-abelian.
    \item[(2)] Suppose $G=\langle X\rangle$ is generated by the set $X$. If $(x,1,\ldots,1)\in \langle F(G,k)\rangle$ for each $x\in X$, then $G$ is non-abelian. 
\end{enumerate}
\end{corollary}

For each $n\geq 3$, we define the \textit{$n$-th dihedral group} as \[D_n=\langle r,a:~r^n=1,a^2=1, ara=r^{-1}\rangle.\] The group $D_n$ has $2n$ elements given by \[D_n=\{1,r,r^2,\ldots,r^{n-1},a,ra,r^2a,\ldots,r^{n-1}a\}.\] Note that $\left(r^ia\right)^2=1$ for any $0\leq i\leq n-1$, and $\left(r^{n/2}\right)^2=1$ whenever $n$ is even. Let $\text{ord}_2(D_n)$ denote the number of elements in $D_n$ of order 2. Then, \[\text{ord}_2(D_n)=\begin{cases}
    n,&\hbox{if $n$ is odd;}\\
    n+1,&\hbox{if $n$ is even.}
\end{cases}\] 

In contrast to Theorem~\ref{the-conf-set-thm}(2), we have the following result.

\begin{theorem}\label{ndihedral-group}
    Let $D_n$ be the dihedral group of order $2n$ with $n\geq 3$. If $n$ is odd, then $F(D_n,2n)$ is a generating set of $D_n^{2n}=D_n\times\cdots\times D_n$ ($2n$ times).
\end{theorem}
\begin{proof}
We have \[F(D_n,2n)=\{\sigma(1,r,r^2,\ldots,r^{n-1},a,ra,r^2a,\ldots,r^{n-1}a):~\sigma\in S_{2n}\}.\] 

For the case $n=3$, we obtain the following equalities: 

\begin{eqnarray*}
(1,r^2,r^2a,ar^2,1,1)&=&(1,r,r^2,a,ra,r^2a)\cdot\\
& & (1,r,a,r^2,ra,r^2a);\\
(r,1,1,1,1,1)&=& (r^2,r^2a,ar^2,1,1,1)^2.
\end{eqnarray*}

Moreover, \begin{eqnarray*}
    (1,1,a,a,1,1)&=& (1,r,r,1,1,1)\cdot\\
    & & (1,r^2,r^2a,ar^2,1,1)\cdot\\
    & & (1,1,1,r,1,1);\\
    (1,r,r^2,a,r,r^2)&=& (1,r,r^2,a,ra,r^2a)\cdot\\
    & & (1,1,1,1,a,a);\\
    (1,1,1,a,1,1)&=& (1,r,r^2,a,r,r^2)\cdot\\
    & &\left(1,r^2,1,1,r^2,1\right)\cdot\\
    & &(1,1,r,1,1,r).
\end{eqnarray*}

Hence, we obtain $(r,1,1,1,1,1), (a,1,1,1,1,1)\in \langle F(D_3,6)\rangle$ and, by Lemma~\ref{lem:ejg}(2),  $F(D_3,6)$ is a generating set of $D_3^6$.
    
Now suppose $n>3$ is odd. The following equalities hold:
$$
\begin{array}{rccccccccl}
1 & r & r^2 & \ldots & r^{n-2} & r^{n-1} & a & ra & \ldots & r^{n-1}a\\
1 & r & r^2 & \ldots & r^{n-2} & a & r^{n-1} & ra & \ldots & r^{n-1}a\\ \hline 
1 & r^2 & (r^2)^2 & \ldots & (r^{n-2})^2 & r^{n-1}a & ar^{n-1} & 1 & \ldots & 1 \\
\end{array}
$$ 

\medskip 
$$
\begin{array}{rccccccccl}
1 & r^2 & (r^2)^2 & \ldots & (r^{n-2})^2 & r^{n-1}a & ar^{n-1} & 1 & \ldots & 1\\
1 & r^2 & (r^{n-2})^2 & \ldots & (r^{2})^2 & r^{n-1}a & ar^{n-1} & 1 & \ldots & 1\\ \hline 
1 & r^4 & 1 &\ldots & 1 & 1 & 1 & 1 &\ldots & 1 \\
\end{array}
$$ 

Here, we use the fact that $n-2-2+1$ is even.  Let $k,j\geq 1$ be integers such that $4j=nk+1$ (recall that $n$ is odd). Observe that \begin{align*}
  (r,1,\ldots,1)&= (r^4,1,\ldots,1)^{j}.  
\end{align*}

Moreover, 

$$
\begin{array}{rccccccccl}
1 & (r^{n-1})^2 & (r^{n-2})^2 & \ldots & (r^{2})^2 & r & 1 & 1 & \ldots & 1\\
1 & r^2 & (r^{2})^2 & \ldots & (r^{n-2})^2 & r^{n-1}a & ar^{n-1} & 1 & \ldots & 1\\
1 & 1 & 1 &\ldots & 1 & 1 & r & 1 & \ldots & 1\\ \hline 
1 & 1 & 1 &\ldots & 1 & a & a & 1 & \ldots& 1 \\
\end{array} 
$$ 

\medskip
$$
\begin{array}{rccccccccccc}
1 & r & r^2 & \ldots & r^{n-2} & r^{n-1} & a & ra & r^2a & \ldots & r^{n-2}a & r^{n-1}a\\
1 & 1 & 1 & \ldots & 1 & 1 & 1 & a & a &\ldots &a & a\\ \hline 
1 & r & r^2 &\ldots & r^{n-2} & r^{n-1} & a & r & r^2 &\ldots & r^{n-2} & r^{n-1} \\
\end{array} 
$$ 

\medskip
$$
\begin{array}{rccccccccccc}
1 & r & r^2 &\ldots & r^{n-2} & r^{n-1} & a & r & r^2 &\ldots & r^{n-2} & r^{n-1}\\
1 & r^{n-1} & r^{n-2} &\ldots & r^2 & r & 1 & r^{n-1} & r^{n-2} & \ldots & r^2 & r\\ \hline 
 1 & 1 & 1 & \ldots & 1 & 1 & a & 1 & 1 &\ldots & 1 & 1\\
\end{array} 
$$

Hence, we obtain $(r,1,\ldots,1), (a,1,\ldots,1)\in \langle F(D_n,2n)\rangle$ and, by Lemma~\ref{lem:ejg}(2),  $F(D_n,2n)$ is a generating set of $D_n^{2n}$.
    
\end{proof}

\begin{remark}\label{rem:even-diedral}
  We use the GAP System to verify whether $\langle F(G,|G|)\rangle=G^{|G|}$ for the dihedral group $D_4$ and the quaternion group $Q_8=(C_2\times C_2)\rtimes(C_2)$, which is a non-abelian group of order $8$ that is not isomorphic to $D_4$. Our computational verifications showed that \[[D_4^{8}:\langle F(D_4,8)\rangle]=4=[Q_8^{8}:\langle F(Q_8,8)\rangle],\] where $[G:H]$ denotes the index of $H$ as a subgroup of $G$.

This shows that the statement \aspas{for any non-abelian finite group $G$, it holds that $\langle F(G,|G|)\rangle=G^{|G|}$} is false. We suspect that $F(D_n,2n)$ is not a generating set of $D_n^{2n}$ for any even $n\geq 4$.
\end{remark}

\section{Applications}\label{sec:applications}

\subsection{Characterization of $F(G,k)$} Let $G$ be a group and let $k\geq 1$. Define the natural projection $$p_{k+1,k}:G^{k+1}\to G^k,~p(g_1,\ldots,g_k,g_{k+1})=(g_1,\ldots,g_k).$$ 
  
  \begin{definition}[Configuration Property]\label{defn:cp}
   Let $C$ be a subset of $G^{k+1}$ (with $k\geq 1$). We say that $C$ has the \textit{configuration property} if it  satisfies the following conditions:
  \begin{itemize}
      \item[(P-1)] $ p_{k+1,k}(C)\cap F(G,k)\neq \emptyset$; and 
      \item[(P-2)] If $ p_{k+1,k}(c)\in F(G,k)$ for some $c\in C$, then $c\in F(G,k+1)$.
  \end{itemize}
  \end{definition}
  
 Any subset of the configuration set $F(G,k+1)$ has the configuration property. In addition, we have the following examples.
  
 \begin{example}\label{exam:subset-X}
\noindent\begin{enumerate}
    \item[(1)] Let $G=\mathbb{Z}_2=\{\overline{0},\overline{1}\}$ and $k=2$. The configuration set is given by \[F(G,k)=\{(\overline{0},\overline{1}),(\overline{1},\overline{0})\}.\] The set $X=\{(\overline{0},\overline{1}),(\overline{1},\overline{1})\}$ is another generating set of $G^2$ that does not have the configuration property, while the subset $\{(\overline{0},\overline{1})\}$ does have the configuration property.
    \item[(2)] For $k\geq 3$, consider the subset \[E_k=\{(g_1,1,\ldots,1),(1,g_2,1,\ldots,1),\ldots,(1,\ldots,1,g_k):~g_i\in G\}\subseteq G^k.\] This is a generating set of $G^k$ and $F(G,k)\cap E_k=\emptyset$. Furthermore, $E_k$ does not have the configuration property (note that $p_{k,k-1}(E_k)=E_{k-1}$).
    \item[(3)] Let $\alpha\in G^k\setminus F(G,k)$ with $|G|\geq k+1$. The set $C(\alpha)=F(G,k)\cup \{\alpha\}$ is a generating set of $G^k$ (by Theorem~\ref{the-conf-set-thm}(1)). Note that $C(\alpha)$ has the configuration property if and only if $p_{k,k-1}(\alpha)\notin F(G,k-1)$.
\end{enumerate}
\end{example} 

Let $(P,\leq)$ be a poset. A \textit{maximum element} is an element that is greater than or equal to all others in $P$. In any poset, there is at most one maximum element. Likewise, a \textit{minimum element} is an element that is less than or equal to all others in $P$, and there is at most one.  

\medskip The following result provides an algebraic-combinatorial characterization of the set $F(G,k)$.

\begin{theorem}[Characterization]\label{charact-k2}  
 Let $G$ be a group. \begin{enumerate}
    \item[(1)] Define the poset $\mathfrak{T}(G)$ whose elements are all subsets $C$ of $G\times G$ that have the configuration property, with $C\leq C'$ if $C\subseteq C'$. Then, $F(G,2)$ is the maximum element of $\mathfrak{T}(G)$.
    \item[(2)] Suppose $k\geq 4$ and $| G| \geq k+1$. Define the poset $\mathfrak{F}(G,k)$ whose elements are all proper subsets $C$ of $G^k$ that satisfy the following conditions:
\begin{itemize}
    \item[(C-1)] The subset $C$ is a generating set of $G^k$ such that $C\cap E_k=\varnothing$;
    \item[(C-2)]  $C\cap S\neq\emptyset$ for every generating set $S$ of $G^k$ having the configuration property. 
\end{itemize}  Again, $C\leq C'$ if $C\subseteq C'$. Then,  $F(G,k)$ is the minimum element of $\mathfrak{F}(G,k)$.
\item[(3)] Suppose $| G| \geq 4$. Define the poset $\mathfrak{H}(G)$ whose elements are all proper subsets $C$ of $G^3$ that satisfy the following conditions:
\begin{itemize}
    \item[(C-1)] The subset $C$ is a generating set of $G^3$ such that $C\cap E_3=\varnothing$;
    \item[(C*-2)]  $C\cap S\neq\emptyset$ for every set $S$ of $G^3$ having the configuration property. 
\end{itemize}  Again, $C\leq C'$ if $C\subseteq C'$. Then,  $F(G,3)$ is the minimum element of $\mathfrak{H}(G)$.
\end{enumerate}  
  \end{theorem}
  \begin{proof}
      \noindent\begin{enumerate}
          \item[(1)] By Definition~\ref{defn:cp}, $F(G,2)\in \mathfrak{T}(G)$. Furthermore, a subset $C\subseteq G\times G$ has the configuration property if and only if $C\subseteq F(G,2)$.
          \item[(2)] By Theorem~\ref{the-conf-set-thm}(1) together with Example~\ref{exam:subset-X}(2), we see that $F(G,k)$ satisfies condition (C-1).

          Let $S\subseteq G^k$ be a generating set with the configuration property. By Definition~\ref{defn:cp}, we have $F(G,k-1)\cap p_{k,k-1}(S)\neq\emptyset$, so there exists an element $x\in S$ such that $p_{k,k-1}(x)\in F(G,k-1)$, and hence $x\in F(G,k)$. This implies  $F(G,k)\cap S\neq\emptyset$, so condition (C-2) is also satisfied. Therefore, $F(G,k)\in \mathfrak{F}(G,k)$.
          
          Now, suppose $C\in \mathfrak{F}(G,k)$. We claim that $F(G,k)\subseteq C$. Let $x\in F(G,k)$ and consider $S_x=E_k\cup\{x\}$. Note that $S_x$ is a generating set of $G^k$ satisfying the configuration property (because $E_{k-1}\cap F(G,k-1)=\varnothing$, here we use the fact that $k\geq 4$). Then, by (C-2), we have that $C\cap S_x\neq\varnothing$. Furthermore, by (C-1), we obtain that $x\in C$. Therefore, $F(G,k)\subseteq C$. 
          \item[(3)] Similarly as item (2), we have that $F(G,3)$ satisfies condition (C-1) and (C*-2). 
          
          Suppose $C\in \mathfrak{h}(G)$. We claim that $F(G,3)\subseteq C$. Let $x\in F(G,3)$ and consider $T_x=L\cup\{x\}$ where $L=\{(1,1,g):~g\in G\}\subseteq E_3$. Note that $T_x$ satisfies the configuration property. Then, by (C*-2), we have that $C\cap T_x\neq\varnothing$. Furthermore, by (C-1), $C\cap L=\varnothing$ and thus we obtain that $x\in C$. Therefore, $F(G,3)\subseteq C$. 
      \end{enumerate}
  \end{proof}

 We now present the following example. 

\begin{example}
Let $G$ be a group. Suppose $k\geq 3$ and $|G|\geq k+1$. Then the set $C(\alpha)=F(G,k)\cup\{\alpha\}$, where $\alpha\notin E_k$ (for instance, $\alpha=(g,g,1,\ldots,1)$ with $g\neq 1$), belongs to the poset $\mathfrak{h}(G)$ (for $k=3$) and $\mathfrak{F}^\ast(G,k)$ (for $k\geq 4$), respectively. 
\end{example}

\subsection{Linear systems in $\mathbb{F}_q$}\label{configuration-group}
Let $p$ be a prime number, and let $\mathbb{F}_q$ be the finite field with $3\leq q=p^n$ elements ($n\geq 1$ integer). Note that $\|\cdot\|\colon \mathbb{F}_q^q\to \mathbb{F}_q$ is a $\mathbb{F}_q$-linear transformation. In this case, by Theorem~\ref{the-conf-set-thm}(2),  \begin{align*}
    \langle F(\mathbb{F}_q,q)\rangle&= \mathrm{Ker}(\|\cdot\|).
\end{align*} The configuration group $\langle F(\mathbb{F}_q,q)\rangle$ is a $\mathbb{F}_q$-vector subspace of $\mathbb{F}_q^{q}$ and \begin{align}\label{dim-conf}
   \dim_{\mathbb{F}_q} \langle F(\mathbb{F}_q,q)\rangle&= q-1.
\end{align}

In particular, the vector spaces (over $\mathbb{F}_q$) $\langle F(\mathbb{F}_q,q)\rangle$ and $\mathbb{F}_q^{q-1}$ are isomorphic. Therefore, the number of vectors in $\langle F(\mathbb{F}_q,q)\rangle$ is $q^{q-1}$.

\medskip We have the following example.

\begin{example}\label{exam:F4}
For \[\mathbb{F}_4=\dfrac{\mathbb{F}_2[x]}{\langle x^2+x+1\rangle},\] we have $\mathbb{F}_4=\{0,1,\alpha,\alpha+1\}$ where $\alpha$ is the equivalence class of $x$. Hence, the configuration set is \[F(\mathbb{F}_4,4)=\{\sigma(0,1,\alpha,\alpha+1):~\sigma\in S_4\},\] and \[\langle F(\mathbb{F}_4,4)\rangle=\mathrm{Ker}(\|\cdot\|).\]
\end{example}

\medskip Consider a homogeneous system of $m$ linear equations of the form  \[a_{j,1}x_1+\cdots+a_{j,k}x_k=0,\] for $j=1,\ldots,m$, with coefficients $a_{j,i}\in\mathbb{F}_q$. Suppose $2\leq m\leq q$, and that \[a_i=(a_{1,i},\ldots,a_{m,i})\in \langle F(\mathbb{F}_q,m)\rangle\] for $i=1,\ldots,k$, with $a_i\neq a_j$ for $i\neq j$ (i.e., $(a_1,\ldots,a_k)\in F\left( \langle F(\mathbb{F}_q,m)\rangle,k\right)$). For the case $m=k=q\geq 3$, the following result provides a method for constructing homogeneous linear systems in $\mathbb{F}_q$ that admit nontrivial solutions. 

\medskip We recall that a homogeneous linear system $Ax=0$ over any field (in this case, over $\mathbb{F}_q$) admits a nontrivial solution in $\mathbb{F}_q^q$ if and only if the matrix $A$ is not invertible over $\mathbb{F}_q$, or equivalently, $\mathrm{det}(A)=0$ in $\mathbb{F}_q$. Furthermore, the number of solutions to $Ax=0$ over $\mathbb{F}_q$ is $q^{q-\mathrm{rank}(A)}$.  

\begin{proposition}\label{nontrivial-solution}
Let $p$ be a prime number, and let $3\leq q=p^n$ ($n\geq 1$ integer). Let $a_1,\ldots,a_q\in \langle F(\mathbb{F}_q,q)\rangle$ with $a_i\neq a_j$ for $i\neq j$ (i.e., $(a_1,\ldots,a_q)\in F\left(\langle F(\mathbb{F}_q,q)\rangle,q\right)$). Let $A=(a_1^T \cdots a_q^T)\in \text{M}_q(\mathbb{F}_q)$ (i.e., $A$ is a $q\times q$ matrix with entries in $\mathbb{F}_q$) where the $i$-th column is $a_i^T$. Then, the homogeneous linear system \[Ax=0\] admits a nontrivial solution in $\mathbb{F}_q^q$, i.e., a solution $x\in \mathbb{F}_q^q$ with $x\neq 0$. Furthermore, the number of such non-invertible matrices $A=(a_1^T \cdots a_q^T)$, i.e., the number of elements in the configuration set $F\left(\langle F(\mathbb{F}_q,q)\rangle,q\right)$, is equal to $(q^{q-1})!/\left(q^{q-1}-q\right)!$.
\end{proposition}
\begin{proof}
Suppose $a_i=(a_{1,i},\ldots,a_{q,i})$ and $x=(x_1,\ldots,x_q)^T$. The equation $Ax=0$ is equivalent to \[
x_1a_1+\cdots+x_qa_q=0.
\] Since $\dim_{\mathbb{F}_q}\langle F(\mathbb{F}_q,q)\rangle=q-1$ (by (\ref{dim-conf})), the vectors $a_1,\ldots,a_q$ are linearly dependent. Therefore, there exists a nonzero solution $x=(x_1,\ldots,x_q)\in \mathbb{F}_q^q$.

Furthermore, since the number of vectors in $\langle F(\mathbb{F}_q,q)\rangle$ is $q^{q-1}$, the number of ordered $q$-tuples of distinct vectors in  $F\left(\langle F(\mathbb{F}_q,q)\rangle,q\right)$ is $(q^{q-1})!/\left(q^{q-1}-q\right)!$ (see (\ref{order-conf})).
\end{proof}

We consider the following example.

\begin{example}
\noindent \begin{enumerate}
    \item[(1)] Let $a_1=(0,1,2), a_2=(1,0,2)$, and $a_3=(0,2,1)$ (note that $(a_1,a_2,a_3)\in F\left(\langle F(\mathbb{F}_3,3)\rangle,3\right)$). The equation \[x_1a_1+x_2a_2+x_3a_3=0\] admits a nontrivial solution given by $x_1=1,x_2=0$, and $x_3=1$. Note that  $(1,0,1)\notin \langle F(\mathbb{F}_3,3)\rangle$.
    \item[(2)] Let $a_1=(0,1,2), a_2=(2,0,1)$, and $a_3=(1,2,0)$. Then the equation \[x_1a_1+x_2a_2+x_3a_3=0\] admits a nontrivial solution given by $x_1=x_2=x_3=1$. In this case, the solution $(1,1,1)\in \langle F(\mathbb{F}_3,3)\rangle$.
    \item[(3)] Let $a_1=(0,1,2), a_2=(1,0,2)$, and $a_3=(0,0,1)$. Then, the equation \[x_1a_1+x_2a_2+x_3a_3=0\] does not  admit a nontrivial solution. Note that $a_3=(0,0,1)\notin \langle F(\mathbb{F}_3,3)\rangle$.
     \item[(4)] Let $a_1=(0,1,2), a_2=(0,2,1)$, and $a_3=(0,0,1)$. Then the equation \[x_1a_1+x_2a_2+x_3a_3=0\] admits a nontrivial solution given by $x_1=1$, $x_2=1$, and $x_3=0$. Again, note that $a_3=(0,0,1)\notin \langle F(\mathbb{F}_3,3)\rangle$.
\end{enumerate}
\end{example}

\subsection{Cayley graphs}\label{cayley-graph}

We recall the definition of a Cayley graph from \cite{hamidoune1991}. Let $G$ be a group and let $S$ be a subset of $G\setminus\{1\}$ such that $S=S^{-1}$. The \textit{Cayley graph} of $G$ with respect to $S$ is the undirected graph \[\mathrm{Cay}(G,S)=(G,E),\] where $E=\{\{x,y\}| ~x^{-1}y\in S\}$. Notice that $x^{-1}y\in S$ if and only if $y^{-1}x\in S$, because $S=S^{-1}$. We write $xy$ for the edge $\{x,y\}$. 

\medskip A sequence of distinct vertices $[x_1,x_2,\ldots,x_k]$ such that $x_ix_{i+1}$ is an edge for each $1\leq i\leq k-1$  is called a \textit{path} in $\mathrm{Cay}(G,S)$ from $x_1$ to $x_k$ of length $k-1$. Of course, we include the case $k=1$, where the path consists of a single vertex.

\medskip Let $G$ be a group and let $S$ be a subset of $G\setminus\{1\}$ such that $S=S^{-1}$. A factorization of an element $x\in G\setminus\{1\}$ with respect to $S$ determines a path in $\mathrm{Cay}(G,S)$ of the same length from $1$ to $x$, and vice versa. Let $x=x_1x_2\cdots x_\ell$, with $x_i\in S$, be a factorization of $x$. Then the path \[[1,x_1,x_1x_2,\ldots,x_1x_2\cdots x_{\ell-1},x_1x_2\cdots x_{\ell-1}x_\ell]\] is from $1$ to $x$ and has length $\ell$. Conversely, suppose $[1,x_1,\ldots,x_\ell]$ is a path from $1$ to $x_\ell=x$ of length $\ell$, then \[x=x_1' x_2'\cdots x_\ell',\] where $x_1'=x_1$, and for $i=2,\ldots,\ell$, we have $x_i'=x_{i-1}^{-1}x_{i}$. Note that each $x_i'\in S$. 

\medskip Hence, $S$ is a generating set of $G$ if and only if the Cayley graph $\mathrm{Cay}(G,S)$ is \textit{path-connected}; that is, for any vertices $a,b\in G$, there exists a path in $\mathrm{Cay}(G,S)$ from $a$ to $b$. Furthermore, the index $[G:\langle S\rangle]$ coincides with the number of path-connected components of $\mathrm{Cay}(G,S)$, that is, \[[G:\langle S\rangle]=|\pi_0\mathrm{Cay}(G,S)|,\] where $\pi_0\mathrm{Cay}(G,S):=G/\sim=\{[g]:~g\in G\}$ is the quotient set with respect to the equivalence relation $\sim$ given by: for $a,b\in G$, we have $a\sim b$ if and only if there exists a path from $a$ to $b$. Indeed, the map $G/\langle S\rangle\to G/\sim$, given by $g+\langle S\rangle\mapsto [g]$,  is a bijection. 

\medskip Theorem~\ref{the-conf-set-thm} implies the following statement.

\begin{theorem}\label{cayley-conf}
Let $G$ be a group. 
\begin{enumerate}
    \item[(1)] If $k=2$ or $| G| \geq k+1\geq 4$, then the Cayley graph $\mathrm{Cay}(G^k,F(G,k))$ is path-connected.
    \item[(2)] If $G$ is finite abelian with $| G| =k\geq 3$, then the Cayley graph $\mathrm{Cay}(G^k,F(G,k))$ is not path-connected. Furthermore, \[\left|\pi_0\mathrm{Cay}(G^k,F(G,k))\right|=\left|G/\left\langle \sum_{g\in G}g\right\rangle\right|.\] For instance, $|\pi_0\mathrm{Cay}(G^k,F(G,k))|= k$ whenever $k$ is odd.
\end{enumerate}
\end{theorem}
\begin{proof}
 \noindent\begin{enumerate}
        \item[(1)] This is a direct consequence of Theorem~\ref{the-conf-set-thm}(1).
        \item[(2)] Since $\left\langle F(G,k)\right\rangle=\|\cdot\|^{-1}\left(\left\langle \sum_{g\in G}g\right\rangle\right)$ (see the proof of Theorem~\ref{the-conf-set-thm}(2)), the map \[G^k/\langle F(G,k)\rangle\to G/\left\langle \sum_{g\in G}g\right\rangle,\] given by $a+\langle F(G,k)\rangle\mapsto ||a||+\langle g_2+\cdots+g_k\rangle$, is an isomorphism. Therefore, we obtain our result.   
    \end{enumerate}
\end{proof}

Let $\mathrm{Cay}(G,S)$ and $\mathrm{Cay}(H,T)$ be Cayley graphs. If $f:G\to H$ is a homomorphism such that $f(S)\subseteq T$, then \[f:\mathrm{Cay}(G,S)\to \mathrm{Cay}(H,T),\] that is, $f$ is a graph homomorphism from $\mathrm{Cay}(G,S)$ to $\mathrm{Cay}(H,T)$. Hence, we have the following observation.

\begin{remark}\label{graph-bet-cay}
Let $G$ be a finite abelian group with $|G|=k\geq 3$ even such that $\sum_{g\in G}g\neq 0$. Suppose that $g_2,\ldots,g_k$ are all the non-identity distinct elements of $G$, then $g_2+\cdots+g_k$ has order 2 (see the proof of Theorem~\ref{the-conf-set-thm}(2)). Then $\{g_2+\cdots+g_k\}=\{g_2+\cdots+g_k\}^{-1}$, so we can consider the Cayley graph $\mathrm{Cay}\left(G, \{g_2+\cdots+g_k\}\right)$. Note that the norm map $||\cdot||:G^k\to G$ is a graph homomomorphism from $\mathrm{Cay}(G^k,F(G,k))$ to $\mathrm{Cay}\left(G, \{g_2+\cdots+g_k\}\right)$ because $||\cdot||:G^k\to G$ is a homomorphism and $||\cdot||\left(F(G,k)\right)=\{g_2+\cdots+g_k\}$.  
\end{remark}


Let $\chi$ denote the usual chromatic number \cite[p. 6]{hell2004}. We obtain the following statement. 

\begin{proposition}
  Let $G$ be a finite abelian group with $|G|=k\geq 3$ such that $\sum_{g\in G}g\neq 0$. Then \[\chi\left(\mathrm{Cay}\left(G^k,F(G,k)\right)\right)=2.\]  
\end{proposition}
\begin{proof}
  From Remark~\ref{graph-bet-cay}, we have a graph homomorphism \[\|\cdot\|\colon \mathrm{Cay}\left(G^k,F(G,k)\right)\to \mathrm{Cay}\left(G, \{a\}\right),\] where $a:=\sum_{g\in G}g$. Hence, by \cite[Corollary 1.8, p. 7]{hell2004}, \[\chi\left(\mathrm{Cay}\left(G^k,F(G,k)\right)\right)\leq \chi\left(\mathrm{Cay}\left(G, \{a\}\right)\right).\] Since $a$ has order $2$, the Cayley graph $\mathrm{Cay}\left(G, \{a\}\right)$ connects each element $x\in G$ to $ax$. Because $a^2=1$, applying $a$ twice returns to the original element:\[x\to ax \to a(ax)=x.\] Thus, every vertex $x$ is paired with exactly one distinct vertex $ax$, and these pairs form disjoint edges. Therefore, \[\chi\left(\mathrm{Cay}\left(G, \{a\}\right)\right)\leq 2.\] On the other hand, by Theorem~\ref{cayley-conf}(2), the number of path-connected components of $\mathrm{Cay}\left(G^k,F(G,k)\right)$ is $k/2$. Hence, \[2\leq \chi\left(\mathrm{Cay}\left(G^k,F(G,k)\right)\right).\] Therefore, \[\chi\left(\mathrm{Cay}\left(G^k,F(G,k)\right)\right)= \chi\left(\mathrm{Cay}\left(G, \{a\}\right)\right)=2.\]
\end{proof}

\subsection{The set $F(G\setminus\{1\},k)$}\label{f-punctured} 

Let $k\geq 2$ be an integer and let $G$ be a group such that $|G|\geq k+1$. By Theorem~\ref{the-conf-set-thm}(1), the configuration set $F(G,k)$ is a generating set of $G^k$. The purpose of this section is to study the configuration set $F(G\setminus\{1\},k)\subsetneq F(G,k)$.

\medskip If $k\geq 2$ is an integer and $G$ is a group such that $|G|\geq k+1$, then the configuration set $F(G,k+1)$ is in bijection with the product \[G\times F(G\setminus\{1\},k),\] via the map $(g_0,g_1,\ldots,g_k)\mapsto (g_0,g_1g_0^{-1},\ldots,g_kg_0^{-1})$, whose inverse is given by $(g_0,g_1,\ldots,g_k)\mapsto (g_0,g_1g_0,\ldots,g_kg_0)$.  

\medskip We consider the surjective map $\varphi:G^{k+1}\to G^{k}$ defined by \begin{equation*}
\varphi(x_0,x_1,\ldots,x_k)=\left(x_1x_0^{-1},\ldots,x_kx_0^{-1}\right).
\end{equation*} 

This map captures the commutativity of the group $G$, as shown in the following remark. 

\begin{remark}\label{rem:varphi-hom}
We note that $\varphi$ is a homomorphism if and only if $G$ is abelian. Indeed, suppose that $\varphi$ is a homomorphism and take any $x,y\in G$. Then: \begin{align*}
 (yx,x,\ldots,x) &=\varphi(x^{-1},y,1,\ldots,1)\\
 &=\varphi\left((x^{-1},1,1,\ldots,1)\cdot (1,y,1,\ldots,1)\right)\\
 &=\varphi(x^{-1},1,1,\ldots,1)\cdot\varphi(1,y,1,\ldots,1)\\
 &=(1\cdot x,1\cdot x,\ldots,1\cdot x)\cdot (y\cdot 1^{-1},1\cdot 1^{-1},\ldots,1\cdot 1^{-1})\\
 &=(xy,x,\ldots,x)\\   
\end{align*} and thus $xy=yx$, so $G$ is abelian. The other implication holds immediately. 
\end{remark}

We now state the following result.

\begin{proposition}\label{prop:graph-homo}
 Let $k\geq 2$ be an integer and let $G$ be an abelian group such that $|G|\geq k+1$. Then the map $\varphi:G^{k+1}\to G^{k}$ is a graph homomorphism from $\mathrm{Cay}\left(G^{k+1},F(G,k+1)\right)$ to $\mathrm{Cay}\left(G^k,F(G\setminus\{0\},k)\right)$. 
\end{proposition}
\begin{proof}
    This follows because $\varphi$ is a homomorphism (by Remark~\ref{rem:varphi-hom}) and $\varphi(F(G,k+1))= F(G\setminus\{0\},k)$.
\end{proof}

We now obtain the following result.

\begin{theorem}\label{thm:punctured-group}
 Let $k\geq 2$ be an integer and let $G$ be an abelian group such that $|G|\geq k+1$. 
\begin{enumerate}
    \item[(1)] If $| G| \geq k+2$, then $F(G\setminus\{0\},k)$ is a generating set of $G^k$. In particular, if $G$ is an infinite abelian group, the subgroup generated by $F(G\setminus\{0\},k)$ equals $G^k$ for any $k\geq 2$.
    \item[(2)] If $|G|= k+1$, then $F(G\setminus\{0\},k)$ is not a generating set of $G^k$. 
\end{enumerate}   
\end{theorem}
\begin{proof}
    \noindent\begin{enumerate}
        \item[(1)] By Theorem~\ref{cayley-conf}(1), the Cayley graph $\mathrm{Cay}\left(G^{k+1},F(G,k+1)\right)$ is path-connected. In addition, by Proposition~\ref{prop:graph-homo}, we have a vertex-surjective graph homomorphism \[\varphi:\mathrm{Cay}\left(G^{k+1},F(G,k+1)\right)\to \mathrm{Cay}\left(G^k,F(G\setminus\{0\},k)\right).\] Hence, $\mathrm{Cay}\left(G^k,F(G\setminus\{0\},k)\right)$ is path-connected, which is  equivalent to saying that $F(G\setminus\{0\},k)$ is a generating set of $G^k$.
        \item[(2)] This follows analogously to the proof of Theorem~\ref{the-conf-set-thm}(2). Indeed, we have \[\langle F(G\setminus\{0\},k)\rangle\subseteq  \|\cdot\|^{-1}\left(\left\langle g_2+\cdots+g_k\right\rangle\right)\] where $\langle g_2+\cdots+g_k\rangle$ is a proper subgroup of order at most 2. 
    \end{enumerate}
\end{proof}

We obtain the following example.

\begin{example}
 Let $k\geq 2$ be an integer, and let $G$ be an abelian group such that $|G|\geq k+1$. \begin{enumerate}
     \item[(1)] Suppose that $|G|\geq k+2$. Then, by Theorem~\ref{thm:punctured-group}, $F(G\setminus\{0\},k)$ is a generating set of $G^k$. Since $F(G\setminus\{0\},k)\cap E_k=\emptyset$ and $F(G\setminus\{0\},k)\subsetneq F(G,k)$, there exists a subset $S$ of $G^k$ with the configuration property such that $F(G\setminus\{0\},k)\cap S=\emptyset$ (see Theorem~\ref{charact-k2}(2),(3)). For $k\geq 4$ such set $S$ is also a generating set.
     \item[(2)] Suppose that $|G|= k+1$. If $F(G\setminus\{0\},k)$ is a generating set of $G^k$, then $G$ must be  non-abelian.
 \end{enumerate} 
\end{example}

Finally, we propose the following directions for future research. 

\begin{remark}[Future Work]\label{rem:future-work}
\noindent\begin{enumerate}
\item[(1)] Based on Remark~\ref{rem:even-diedral}, we state the following conjecture: The configuration set $F(D_n,2n)$ is not a generating set of $D_n^{2n}$ for any even $n\geq 4$. 
     \item[(2)] Let $k\geq 2$ be an integer and let $G$ be a finite abelian group such that $|G|= k+1$. In view of  Theorem~\ref{thm:punctured-group}(2), we propose to study the subgroup generated by $F(G\setminus\{0\},k)$. 
\end{enumerate}   
\end{remark}

\section*{Conclusion}
We have studied the configuration set $F(G,k)$ for a general (abstract) group $G$ and integer $k\geq 2$, viewed as a subset of $G^k$, and explored various algebraic and combinatorial properties of this set. 

\section*{Conflict of Interest Statement}
On behalf of all authors, the corresponding author states that there is no conflict of interest.

\bibliographystyle{plain}

\end{document}